\documentclass[a4paper,11pt]{article}
\usepackage{amsmath,amscd,amssymb}
\usepackage[frenchb]{babel}
\usepackage{graphicx,psfrag}
\usepackage{xypic}
\usepackage{pstricks}
\usepackage{pst-plot}
\input xy
\xyoption{all}

\def \tha #1#2{\noi{\bf#1{\uppercase{\footnotesize{#2}}}}}

\newtheorem{theor}{\tha{T}{h{\'e}or{\`e}me}}[section]
\newenvironment{theo}{
  \begin{theor}\hs -0.2 cm {\bf .} ---  }
{  \end{theor}}

\newtheorem{propo}[theor]{\tha{P}{roposition}}
\newenvironment{prop}{
  \begin{propo}\hs -0.2 cm {\bf .} ---  }
{  \end{propo}}

\newtheorem{lemma}[theor]{\tha{L}{emme}}
\newenvironment{lemm}{
  \begin{lemma}\hs -0.2 cm {\bf .} ---  }
{  \end{lemma}}

\newtheorem{fait}[theor]{\tha{F}{ait}}
\newenvironment{fact}{
  \begin{fait}\hs -0.2 cm {\bf .} ---  }
{  \end{fait}}

\newtheorem{defini}[theor]{\tha{D}{{\'e}finition}}

\newtheorem{corollaire}[theor]{\tha{C}{orollaire}}
\newenvironment{coro}{
  \begin{corollaire}\hs -0.2 cm {\bf .} ---  }
{  \end{corollaire}}

\newtheorem{exemple}[theor]{\sc{Exemple}}

\newtheorem{remarq}[theor]{\sc{Remarque}}
\newenvironment{rema}{
  \begin{remarq}\hs -0.2 cm {\bf .} ---  }
{  \end{remarq}}

\newenvironment{preu}{\noi\sc{Preuve}\hs 0.1 cm 
--- \rm }{\hfill$\Box$\vs 0.2 cm}

\def \vs {\vskip}
\def \hs {\hskip}
\def \noi {\noindent}

\def \cad {c'est-{\`a}-dire }

\def \oo {{\cal O}}
\def \L {{\cal L}}

\def \a {{\alpha}}
\def \b {{\beta}}
\def \ga {{\gamma}}

\def \vp {\varphi}

\def \C {{\mathbb C}}

\def \p {{\mathbb P}}

\def \G {{\mathbb G}}

\def \fll {\longrightarrow}

\def \Hom {{\rm Hom}}
\def \pr {{\rm pr}}

\def \pic {{\rm Pic}}

\def \MorC #1#2{{\bf{Hom}}_{#1}(C,#2)}

\def \sca #1#2{\left\langle#1,#2\right\rangle}
\def \scal #1#2{\langle #1,#2 \rangle}

\def \Xt {{\widetilde{X}}}

\def \ft {{\widetilde{f}}}
\def \gt {{\widetilde{g}}}

\def \at {{\widetilde{\a}}}
\def \wt {{\widetilde{w}}}

\def \pt {{\widetilde{\p}}}

\def \cw {\overline{w}}
\def \cs {\overline{s}}

\def \app {A_1^+(\Xt_\wp)}

\def \wp {{F_{\bullet}}}

\begin{document}

~
\vs -1 cm

\centerline{\large{\uppercase{\bf{Courbes elliptiques}}}}
\centerline{\large{\uppercase{\bf{sur la vari{\'e}t{\'e}
        spinorielle}}}}

\vs 0.5 cm

\centerline{\Large{Nicolas \textsc{Perrin}}}

\vs 1 cm

\centerline{\Large{\bf Introduction}} 

\vs 0.3 cm

Dans cet article, nous d{\'e}montrons l'irr{\'e}ductibilit{\'e} du sch{\'e}ma des
morphismes d'une courbe elliptique vers la vari{\'e}t{\'e}
spinorielle. 

\vs 0.2 cm

Soit $G$ le
groupe $SO(2n)$ et $P$ le sous-groupe parabolique maximal de $G$
associ{\'e} --- avec les notation de N. Bourbaki \cite{bourb} --- {\`a} la
racine simple $\a_n$. La vari{\'e}t{\'e} homog{\`e}ne $G/P$ est une composante
connexe de la grassmannienne
des sous-espaces totalement isotropes maximaux de dimension $n$ dans
$\C^{2n}$ muni d'une forme quadratique non d{\'e}g{\'e}n{\'e}r{\'e}e. 
Nous noterons $X$ cette vari{\'e}t{\'e} qui est lisse de dimension
$N=\displaystyle{\frac{n(n-1)}{2}}$.

\vs 0.2 cm

Soit $\a\in A_1(X)$ une classe de 1-cycles sur $X$ et $\MorC{\a}{X}$
le sch{\'e}ma des morphismes de classe $\a$ d'une courbe lisse $C$ vers
$X$. Nous montrons le r{\'e}sultat suivant :

\begin{theo}
\label{main}
  Soit $C$ une courbe elliptique lisse et soit $\a\in A_1(X)$ une
  classe de 1-cycles sur $X$ de degr{\'e} $d$. D{\`e}s que $d\geq n-1$, le
  sch{\'e}ma $\MorC{\a}{X}$ 
est irr{\'e}ductible de la dimension attendue
$$\int_\a c_1(X)=2(n-1)d.$$
\end{theo}

Outre l'int{\'e}r{\^e}t propre de ce r{\'e}sultat, notre {\'e}tude est motiv{\'e}e par
une question d'Atanas Iliev et Dimitri Markushevich. Dans leur
article \cite{Marku}, ils montrent que l'espace des modules
$M_{X_{12}}(2,1,6)$ des fibr{\'e}s vectoriels de rang 2 et classe de Chern
$(1,6)$ sur le volume de Fano $X_{12}$ d'indice 1 et de degr{\'e} 12
s'identifie par la construction de Serre aux courbes elliptiques de
degr{\'e} 6 sur $X_{12}$. Cependant $X_{12}$ est obtenue comme section
lin{\'e}aire de la vari{\'e}t{\'e} spinorielle de dimension 10. Cette description
de $M_{X_{12}}(2,1,6)$, notre r{\'e}sultat pour $(n,d)=(5,6)$ et un
argument de monodromie leur permet de montrer l'irr{\'e}ductibilit{\'e} de
$M_{X_{12}}(2,1,6)$.

\vs 0.2 cm

Pour montrer le th{\`e}or{\`e}me \ref{main}, nous utilisons la r{\'e}solution de
Bott-Samelson $\pi:\Xt_\wp\to X$ qui est un morphisme propre et
birationnel d{\'e}fini pour tout drapeau complet $\wp$, voir paragraphe
\ref{Bott}. Le principe est de relever un morphisme
$f:C\to X$ en un morphisme $\ft:C\to\Xt$. On {\'e}tudie alors le sch{\'e}ma
$\MorC{\at}{\Xt_\wp}$ pour les classes $\at$ telle que
$\pi_*\at=\a$. La vari{\'e}t{\'e} $\Xt_\wp$ peut {\^e}tre vue comme une tour de
fibrations en droites projectives (cf. paragraphe \ref{Bott}) ce qui
permet de raisonner par r{\'e}currence en se ramenant {\`a} l'{\'e}tude d'une
telle fibration.

Cependant, {\`a} la diff{\'e}rence du cas des courbes rationnelles que nous
avons trait{\'e} dans \cite{perrin} ou \cite{Perrin}, nous avons besoin
d'imposer des conditions de positivit{\'e} plus fortes --- not{\'e}es
$(\star)$ --- sur les classes de 1-cycles $\at\in A_1(\Xt)$ 
pour obtenir l'irr{\'e}ductibilit{\'e} du sch{\'e}ma $\MorC{\at}{\Xt_\wp}$,
voir paragraphe \ref{irredu}. Nous montrons aux propositions
\ref{intersection} et \ref{reste} que pour (essentiellement) tout
morphisme $f\in\MorC{\a}{X}$, il existe un drapeau complet $\wp$ tel
que $f$ se rel{\`e}ve dans $\Xt_\wp$ en un morphisme $\ft$ satisfaisant
les conditions $(\star)$. On conclue alors {\`a} l'irr{\'e}ductibilit{\'e} en faisant
varier $\wp$.

\begin{rema}
(\i) A. Brugui{\`e}res \cite{Bruguieres} a montr{\'e} un r{\'e}sultat d'irr{\'e}ductibilit{\'e}
semblable au th{\'e}or{\`e}me \ref{main} pour les courbes elliptiques trac{\'e}es
sur les grassmanniennes (voir {\'e}galement \cite{PEell}).
Le cas des courbes trac{\'e}es sur des quadriques est trait{\'e} par
E. Ballico dans \cite{Ballico}.

(\i\i) La technique pr{\'e}sent{\'e}e ici a l'avantage de pouvoir se g{\'e}n{\'e}raliser {\`a}
toutes les vari{\'e}t{\'e}s homog{\`e}nes minuscules et permettrait de retrouver
les r{\'e}sultats du th{\'e}or{\`e}me \ref{main}, de A. Brugui{\`e}res et de
E. Ballico de mani{\`e}re unifi{\'e}e. Ceci fait l'objet d'une partie du texte
\cite{PEcarquois}  en pr{\'e}paration.

(\i\i\i) Il semble possible qu'en utilisant de mani{\`e}re plus subtile le
carquois $Q$ associ{\'e} {\`a} la r{\'e}solution de Bott-Samelson (cf. paragraphe
\ref{choix}) le r{\'e}sultat reste vrai sans condition sur le degr{\'e}. 
\end{rema}

\section{R{\'e}solution de Bott-Samelson}
\label{Bott}

Dans ce paragraphe, nous reprenons les notations de 
\cite{small}. Soit $W$ le groupe de Weyl de $G$ et $W_P$ le
sous-groupe de $W$ stabilisant $P$. On note $w_0$ l'{\'e}l{\'e}ment de
longueur maximale dans $W$ et $\cw_0$ sa classe dans $W/W_P$. Il est
bien connu (cf. \cite{stem1}) que $\cw_0$ a une unique {\'e}criture
r{\'e}duite {\`a} relation de commutation pr{\`e}s. Nous la noterons $\wt_0$ et la
fixons sous la forme
suivante :
$$\cw_0=\cs_{\b_1}\cdots\cs_{\b_{N}}$$
o{\`u} les $\b_i$ sont les racines simples d{\'e}finies de la fa{\c c}on
  suivante. On prend les notations de \cite{bourb} pour les racines
  simples qui sont not{\'e}es $\a_1,\cdots,\a_n$. On d{\'e}finit pour tout
  $k\in[1,n]$ l'entier
$$a_k=\frac{(k-1)(2n-k)}{2}$$
et on d{\'e}coupe l'intervalle $[1;N]$ selon les intervalles disjoints
  $[a_k+1;a_{k+1}]$ pour $k\in[1;n-1]$
(remarquons que $a_{k+1}=a_k+n-k$). Alors si $i=a_k+j$ avec
$j\in[1,n-k]$, on pose 
$$\b_i=\left\{
  \begin{array}{l}
\a_{n-j}\ \textrm{ si } j\geq2\\
\a_{n}\ \textrm{ si } j=1 \textrm{ et } n-k \textrm{ est impair}\\
\a_{n-1}\ \textrm{ si } j=1 \textrm{ et } n-k \textrm{ est pair.}\\
  \end{array}\right.$$
{\`A} une telle {\'e}criture r{\'e}duite, nous avons associ{\'e} dans \cite{small} un
carquois $Q$. Nous donnons une repr{\'e}sentation de $Q$ en appendice
ce qui permet de mieux comprendre le choix des notations ci-dessus.

Une fois cette d{\'e}composition fix{\'e}e, {\`a} tout drapeau complet $\wp$, on
peut associer (cf. \cite{Demazure}) une
vari{\'e}t{\'e} de Bott-Samelson $\Xt_\wp$. Cette vari{\'e}t{\'e} peut {\^e}tre
d{\'e}crite par une suite de fibrations
$$\Xt_\wp=X_{N}\stackrel{f_{N}}{\fll}\cdots
\stackrel{f_{2}}{\fll}X_1\stackrel{f_{1}}{\fll}X_0\simeq{\rm
  Spec}(k)$$ 
o{\`u} les $f_i$ sont des fibrations en droite projectives. Chacune de ces
fibrations $f_i$ est munie d'une section $\sigma_i$. Nous noterons
$\xi_i$ la classe du diviseur $\sigma_i(X_{i-1})$ de $X_i$ et par
$T_i$ le fibr{\'e} tangent relatif de la fibration $f_i$. 
Par abus de notations, nous noterons encore $\xi_i$ et $[T_i]$ les
classes des images r{\'e}ciproques de $\sigma_i(X_{i-1})$ et $T_i$ dans
dans $\Xt_\wp$.

\section{Fibrations}

\subsection{Irr{\'e}ductibilit{\'e}}
\label{irredu}

Nous montrons une proposition qui permet de remonter
l'irr{\'e}ductibilit{\'e} du sch{\'e}ma des morphismes {\`a} travers les fibrations en
droites projectives. Cependant, contrairement au cas des courbes
rationnelles (cf. \cite{perrin} prop. 4), la seule condition d'avoir
un degr{\'e} relatif positif ne suffit plus. Soit $C$ une courbe
elliptique lisse.

\begin{prop}
\label{irred}
  Soit $\vp:X\to Y$ une fibration en droites projectives munie d'une
  section $\sigma$ et soit $\at\in A_1(X)$ une classe de
  1-cycles. Notons $T$ le fibr{\'e} tangent relatif et $\xi$ le diviseur
$\sigma(Y)$. Supposons que $\at$ v{\'e}rifie $\at\cdot\xi\geq0$ et
  $\at\cdot(T-\xi)>0$.

Si $\MorC{\vp_*\at}{Y}$ est irr{\'e}ductible, alors $\MorC{\at}{X}$ l'est
{\'e}galement et on a l'{\'e}galit{\'e}  
$$\dim(\MorC{\at}{X})=\dim(\MorC{\vp_*\at}{Y})+\at\cdot T.$$
\end{prop}

\begin{preu}
Notons $E$ un fibr{\'e} vectoriel de rang 2 sur $Y$ tel que
$X=\p_Y(E)$. La section $\sigma$ est donn{\'e}e par une surjection
$E\to L$ o{\`u} $L$ est inversible. Notons $N$ le fibr{\'e} inversible noyau
de cette surjection. 

Nous {\'e}tudions la fibre du morphisme $\MorC{\at}{X}\to\MorC{\vp_*\at}{Y}$
au dessus d'une fl{\`e}che $f:C\to Y$. Un {\'e}l{\'e}ment de cette fibre est donn{\'e}
par un rel{\`e}vement de $f$, c'est-{\`a}-dire par une surjection $f^*E\to M$
o{\`u} $M$ est inversible sur $C$ avec $2\deg(M)-\deg(f^*E)=\at\cdot
T$. Un {\'e}l{\'e}ment de cette fibre est donc donn{\'e} par un fibr{\'e} inversible $M$
de degr{\'e} $d=\frac{\deg(f^*E)+\at\cdot T}{2}$ et par un {\'e}l{\'e}ment
surjectif de $\p\Hom(f^*E,M)$. On a $\at\cdot\xi=\deg(M)-\deg(f^*N)$
et $\at\cdot(T-\xi)=\deg(M)-\deg(f^*L)$. On distingue deux
cas.

Si $\at\cdot\xi>0$ alors $\Hom(f^*E,M)$ est isomorphe {\`a} 
$\Hom(f^*N,M)\oplus\Hom(f^*L,M)$ et est de dimension contante
(par rapport {\`a} $f$) {\'e}gale {\`a} $\at\cdot T$. Le choix de $M$ est
libre. La fibre est donn{\'e}e par le choix de $M$ puis d'une surjection
$f^*E\to M$ \cad par un ouvert (donn{\'e} par la condition de
surjectivit{\'e}) non vide de $\p\Hom(f^*E,M)\times\pic_{d}(C)$. On a donc
une fibration lisse de dimension $\at\cdot T$ au-dessus de
$\MorC{\vp_*\at}{Y}$ d'o{\`u} le r{\'e}sultat.

Si $\at\cdot\xi=0$ alors si $M\not\simeq f^*N$ on a $\Hom(f^*N,M)=0$
donc toute fl{\`e}che $f^*E\to M$ se factorise par $f^*L$ et on ne peut
avoir de fl{\`e}che surjective car
$\deg(M)-\deg(f^*L)=\at\cdot(T-\xi)>0$. Pour tout {\'e}l{\'e}ment de la fibre,
on doit donc avoir un isomorphisme $M\simeq f^*N$. Mais alors comme
$\deg(M)-\deg(f^*L)=\at\cdot(T-\xi)>0$, le faisceau $f^*E$ est
isomorphe {\`a} $M\oplus f^*L$. On a alors un isomorphisme
$\Hom(f^*E,M)\simeq\Hom(f^*N,M)\oplus\Hom(f^*L,M)$. La dimension de
$\Hom(f^*E,M)$ est donc contante (par rapport {\`a} $f$) {\'e}gale
{\`a} $\at\cdot T+1$. La fibre est donn{\'e}e par un ouvert non vide
de $\p\Hom(f^*E,M)$. On a donc encore une fibration lisse de dimension
$\at\cdot T$ au-dessus de $\MorC{\vp_*\at}{Y}$.
\end{preu}

Nous allons appliquer ce r{\'e}sulat {\`a} la vari{\'e}t{\'e} $\Xt_\wp$. Pour cel{\`a}
nous avons besoin de montrer des r{\'e}sultats de positivit{\'e} sur les
classes $\at$ de 1-cycles.

\begin{rema}
C'est un calcul classique d'exprimer les classes $[T_i]$ en fonction
des classes $\xi_i$ (cf. par exemple \cite{Perrin}  prop. 2.11). On a
$$T_i=\sum_{k=1}^i\sca{\ga_k^\vee}{\ga_i}\xi_k$$
o{\`u} les $\ga_i$ sont les racines positives d{\'e}finies par
$\ga_i=s_{\b_1}\cdots s_{\b_{i-1}}(\b_i)$.
\end{rema}

Nous avons montr{\'e} dans \cite{perrin} que pour une vari{\'e}t{\'e} miniscule
(c'est le cas de $X$), on a toujours
$\sca{\ga_k^\vee}{\ga_i}\in\{0;1\}$. Nous calculons certaines de ces
valeurs.

\begin{lemm}
(\i) Supposons que $i=a_k+j$ pour $k\in[2;n-1]$ et $j\in[1,n-k]$, alors
  on a $\scal{\ga_{k-1}^\vee}{\ga_i}=1$. 

(\i\i) Par ailleurs pour tout $i$ et $j$ distincts dans $[1;n-1]$, on
a $\scal{\ga_{j}^\vee}{\ga_i}=1$.
\end{lemm}

\begin{preu}
  (\i) On commence par remarquer que l'on a l'{\'e}galit{\'e} suivante:
$$\scal{\ga_{k-1}^\vee}{\ga_i}=\scal{\b_k^\vee}{s_{\b_k}\cdots
  s_{\b_{i-1}}(\b_i)}.$$
Par ailleurs, un calcul simple donne 
$$s_{\b_{k-1}}\cdots_{\b_{i-1}(\b_i)}= \sum_{u=n-k-j+1}^{n-k}\a_u+
2\sum_{u=n-k+1}^{n-2}\a_u+\a_{n-1}+\a_n$$
o{\`u} la seconde somme est vide lorsque $k=2$. Le r{\'e}sultat en d{\'e}coule.

(\i\i) Dans ce cas, on a $\ga_i=\b_1+\cdots+\b_i$ et
$\ga_j=\b_1+\cdots+\b_j$ ce qui donne le r{\'e}sultat.
\end{preu}

Notons $\app$ l'ensembles des classes $\at\in A_1(\Xt_\wp)$ telles que
$\at\cdot\xi_i\geq0$ pour tout $i\in[1;r]$ et $\at\cdot\xi_i>0$ pour
tout $i\in[1;n-2]$.

\begin{coro}
\label{def-app}
(\i) Si $\at\in\app$, alors on a $(T_i-\xi_i)\cdot\at>0$ pour tout
  $i\in[1;r]$.

(\i\i) En particulier, si $\at\in\app$, alors le sch{\'e}ma
$\MorC{\at}{\Xt_\wp}$ est irr{\'e}ductible.
\end{coro}

\begin{preu}
  (\i) On a
$$(T_i-\xi_i)\cdot\at=\xi_i\cdot\at+ \sum_{k=1}^{i-1}
\sca{\ga_k^\vee}{\ga_i} \xi_k\cdot\at.$$
Comme on a $\xi_k\cdot\at\geq0$ et $\sca{\ga_k^\vee}{\ga_i}\geq0$ pour
tout $k$, on en d{\'e}duit que pour tout $i\in[1,r]$ et tout
$k\in[1,i-1]$, on a $(T_i-\xi_i)\cdot\at\geq\xi_i\cdot\at$ et
$(T_i-\xi_i)\cdot\at\geq\sca{\ga_k^\vee}{\ga_i}\xi_k\cdot\at$. 

Soit $i\in[1,r]$ que l'on {\'e}crit sous la forme $i=a_k+j$ avec
$k\in[1,n-1]$. Si on a $k\geq2$, alors $\sca{\ga_{k-1}^\vee}{\ga_i}=1$
et $(T_i-\xi_i)\cdot\at\geq\xi_{k-1}\cdot\at>0$. Sinon, on a
$\sca{\ga_1^\vee}{\ga_i}=1$ et
$(T_i-\xi_i)\cdot\at\geq\xi_1\cdot\at>0$.

(\i\i) On applique la proposition \ref{irred} {\`a} toutes les
fibrations $f_i$ et on conclue par r{\'e}ccurence.
\end{preu}

Nous allons montrer au paragraphe \ref{choix} que pour un {\'e}l{\'e}ment $f$
g{\'e}n{\'e}ral de $\MorC{\a}{X}$, on peut choisir un drapeau $\wp$ tel que $f$
se rel{\`e}ve en $\ft$ dans $\MorC{\at}{\Xt_\wp}$ avec $\at\in\app$.

\subsection{Dimension}
\label{dimension}

Nous calculons dans ce paragraphe la dimension du sch{\'e}ma
$\MorC{\at}{\Xt\wp}$ pour une classe $\at\in\app$. Ce sch{\'e}ma est
toujours irr{\'e}ductible comme on l'a vu au corollaire \ref{def-app}. La
dimension de ce sch{\'e}ma est 
$$\sum_{i=1}^rT_i\cdot\at=-K_{\Xt_\wp}\cdot\at.$$
Par ailleurs nous avons montr{\'e} dans \cite{small} que le diviseur
canonique s'exprime simplement en termes des diviseurs $\xi_i$ par la
formule suivante :
$$K_{\Xt_\wp}=\sum_{i=1}^r(h(i)+1)\xi_i$$
o{\`u} $h(i)$ est la hauteur du sommet $i$ c'est-{\`a}-dire la longueur du
plus long chemin du sommet $i$ au sommet $r$ (ainsi par exemple
$h(r)=1$, $h(r-1)=2$, $h(1)=2n-3$, etc.). Par ailleurs, nous avons
montr{\'e}, toujours dans \cite{small}, que le diviseur ample $\L$
g{\'e}n{\'e}rateur du groupe de Picard de $X$ se rel{\`e}ve dans $\Xt_\wp$ en : 
$$\pi^*\L=\sum_{i=1}^r\xi_i$$
ce qui donne :
$$-K_{\Xt_\wp}=(h(1)+1)\pi^*\L-\sum_{i=1}^r(h(1)-h(i))\xi_i
=2(n-1)\pi^*\L-\sum_{i=1}^r(h(1)-h(i))\xi_i.$$
Ainsi pour toute classe $\at\in\app$ telle que $\pi_*\at=\a$ o{\`u} $\a$
est de degr{\'e} $d$, on a
$$-K_{\Xt_\wp}\cdot\at = 2(n-1)d-\sum_{i=1}^r(h(1)-h(i))\xi_i\cdot\at\leq
2(n-1)d-\sum_{i=1}^{n-2}(h(1)-h(i)).$$
Et comme pour tout $i\in[1,n-1]$ on a $h(i)=2(n-1)-i$ (cf. le
carquois donn{\'e} en appendice), on obtient l'in{\'e}galit{\'e}
$$\dim\MorC{\at}{\Xt_\wp}=-K_{\Xt_\wp}\cdot\at
\leq2(n-1)d-\frac{(n-2)(n-3)}{2}$$
avec {\'e}galit{\'e} si et seulement si on a $\at\cdot\xi_i=1$ pour tout
$i\in[2,n-2]$ et $\at\cdot\xi_i=0$ pour tout $i\in[n,r]$.
Nous noterons $\at_0$ l'unique classe telle que $\pi_*\at=\a$
v{\'e}rifiant cette condition.

\section{Choix du drapeau $\wp$}
\label{choix}

\subsection{Retour sur Bott-Samelson}

Rappelons qu'{\`a} partir d'une d{\'e}composition r{\'e}duite $\wt_0$ de $\cw_0$
et d'un drapeau complet $\wp$ on peut construire la r{\'e}solution de
Bott-Samelson $\pi:\Xt_\wp\to X$. D'autre part, P. Magyar
\cite{Magyar} a d{\'e}crit la r{\'e}solution de Bott-Samelson comme une
vari{\'e}t{\'e} de configurations ce que nous avons r{\'e}interpr{\'e}t{\'e} dans
\cite{small} gr{\^a}ce au carquois $Q$. 

Dans le cas pr{\'e}sent, cette vari{\'e}t{\'e} de configurations est donn{\'e}e de la
mani{\`e}re suivante. Pour un {\'e}l{\'e}ment $i\in[1,r]$, la racine simple $\b_i$
d{\'e}finit une grassmannienne de sous-espaces totalement isotropes que
nous noterons $\G_{iso}(\b_i,2n)$. Nous noterons $\dim\b_i$ la
dimension des espaces totalement isotropes {\'e}l{\'e}ments de cette
grassmannienne. Pour chaque sommet $i\in[1,r]$ du carquois, on choisit
un {\'e}l{\'e}ment $x_i$ dans $\G_{iso}(\b_i,2n)$. On a donc une famille :
$$(x_i)_{i\in[1,r]}\in\prod_{i=1}^r\G_{iso}(\b_i,2n).$$
Nous imposons maintenant les conditions suivantes sur cette
famille. Soit $i\in[1,r]$ et notons $\delta_1$, $\delta_2$ et
$\delta_3$ les racines simples adjacentes {\`a} $\b_i$ (s'il y a seulement
une seule resp. deux racines simples, on garde $\delta_1$
resp. $\delta_1$ et $\delta_2$). Lorsqu'ils existent, notons $j_1$,
$j_2$ et $j_3$ les sommets du carquois munis d'une fl{\`e}che vers le
sommet $i$ et tels que $\b_{j_k}=\delta_k$ pour $k\in\{1,2,3\}$. On
d{\'e}finit les conditions : 
$$(*)_i\ :\ \textrm{Pour tout } k\in\{1,2,3\}, \textrm{ on a }
\left\{\begin{array}{c} x_i\subset x_{j_k} \textrm{ si }
    \dim\b_i<\dim\delta_k \\ 
x_{j_k}\subset x_{i} \textrm{ si } \dim\b_i>\dim\delta_k. 
\end{array}\right.$$
Si le sommet $j_k$ n'existe pas, alors on remplace dans les conditions
ci-dessus l'{\'e}l{\'e}ment $x_{j_k}$ par l'unique {\'e}l{\'e}ment du drapeau complet
$\wp$ appartenant {\`a} la grassmanienne $\G_{iso}(\delta_k,2n)$. 

La vari{\'e}t{\'e} de Bott-Samelson s'interpr{\`e}te alors en terme de vari{\'e}t{\'e}s de
configurations par
$$\Xt_\wp=\left\{(x_i)_{i\in[1,r]}\in\prod_{i=1}^r\G_{iso}(\b_i,2n)\
  /\ (x_i)_{i\in[1,r]} \textrm{v{\'e}rifie les conditions } (*)_i  \textrm{
  pour tout } i\in[1,r] \right\}.$$
Le morphisme de $\pi:\Xt_\wp\to X$ est donn{\'e} par la projection
$(x_i)_{i\in[1,r]}\mapsto x_r$. Comme le
  morphisme $\pi$ est birationnel, pour un {\'e}l{\'e}ment $x_r$ g{\'e}n{\'e}ral dans
  $X$, on peut exprimer tous les {\'e}l{\'e}ments de la famille
  $(x_i)_{i\in[1,r]}$ en fonction de $x_r$ et du drapeaux $\wp$. 
Notons $V$ l'unique {\'e}l{\'e}ment du drapeau complet $\wp$ contenu dans
$\G_{iso}(\b_n,2n)$, on a $V=F_{n-1}$ si $n$ est pair et $V=W_n$
sinon. On a alors par exemple
$$x_{n-1}=x_r\cap V.$$
L'{\'e}l{\'e}ment $x_{n-1}$ est une droite vectorielle de $V$ \cad un point de
$\p(V)$.
On d{\'e}duit de $x_{n-1}$ les points $x_i$ pour $i\in[1,n-2]$ par la formule:
$$x_i=x_{n-1}+F_{i-1}.$$
Ce sont des sous-espaces vectoriels de dimension $i$ de $V$ \cad des
{\'e}l{\'e}ment de la grassmannienne $\G(i,V)$.

\subsection{Une projection de $X$}

Il est clair que l'{\'e}l{\'e}ment $x_{n-1}$ est bien d{\'e}fini {\`a} partir de $x_r$
d{\`e}s que $x_r$ rencontre $V$ en dimension
exactement 1. L'ensemble $U$ des points de $X$ o{\`u} c'est le cas est un
ouvert (une orbite sous le stabilisateur de $V$)
dont le compl{\'e}mentaire est de codimension 3. Sur l'ouvert $U$, on a
donc un morphisme $p:U\to\p(V)$ d{\'e}fini par $p(x_r)=x_{r}\cap V$.
Ce morphisme est un cas particulier de ceux {\'e}tudi{\'e}s dans \cite{perrin}
ou \cite{PEell}. On a ainsi le fait suivant :

\begin{fact}
  Le morphisme $p:U\to\p(V)$ permet de r{\'e}aliser $U$ comme l'espace
  total du fibr{\'e} vectoriel $\Lambda^2\left(T_{\p(V)}(-1)\right)$.
\end{fact}

Nous consid{\'e}rons par ailleurs la vari{\'e}t{\'e} image de $\Xt_\wp$ par la
projection $q$ sur les $n-1$ premiers termes, \cad l'image de la
projection 
$$q:\Xt_\wp\to\prod_{i=1}^{n-1}\G_{iso}(\b_i,2n).$$
Cette image est la vari{\'e}t{\'e} $X_{n-1}$ d{\'e}finie au paragraphe
\ref{Bott}. Nous la noterons $\pt_\wp(V)$ car c'est la
r{\'e}solution de Bott-Samelson $\pi':\pt_\wp(V)\to\p(V)$ de l'espace
projectif $\p(V)$ muni du drapeau complet obtenu comme la trace de
$\wp$ sur $V$. Nous avons alors le diagramme suivant
$$\xymatrix{&\Xt_\wp\ar[d]^{\pi}\ar[r]^{q}&\pt_\wp(V)\ar[d]^{\pi'}\\
U_0\ar@{^{(}->}[r]\ar@{^{(}->}[rd]\ar[ur]&X &\p(V)\\
&U\ar@{^{(}->}[u]^{i}\ar[ur]^{p}&}$$
o{\`u} $U_0$ est la cellule de Schubert associ{\'e}e au drapeau complet
$\wp$. Au dessus de $U_0$ l'application $\pi$ est un isomorphisme.

\begin{rema}
  La vari{\'e}t{\'e} $\pt_\wp(V)$ est munie en tant que vari{\'e}t{\'e} de
  Bott-Samelson d'une base des classes de diviseurs qui est donn{\'e}e par
  l'image des classes $\xi_i$ pour $i\in[1,n-1]$. Par abus de
  notations nous noterons encore $\xi_i$ ces classes de diviseurs dans
  $\pt_\wp(V)$.
\end{rema}

\subsection{Choix du drapeau}

Fixons un morphisme $f:C\to X$. Nous commen{\c c}ons par fixer l'espace $V$
du drapeau, \cad l'espace $F_{n-1}$ si $n$ est pair et l'espace $F_n$
dans le cas contraire. Nous choisissons cet espace de sorte que la
courbe $f(C)$
soit contenue dans $U$. Ceci est possible --- c'est une application du
th{\'e}or{\`e}me de Bertini prouv{\'e} par S. Kleiman \cite{kleiman} --- 
gr{\^a}ce {\`a} l'action du groupe et
le fait que le compl{\'e}mentaire de $U$ est de codimension 3.
Un espace $V$ g{\'e}n{\'e}ral convient.

Nous montrons maintenant la proposition suivante :

\begin{prop}
\label{intersection}
  Si la courbe $p\circ f(C)$ n'est pas contenue dans un espace
  lin{\'e}aire de codimension 2, alors on peut compl{\'e}ter
  $V$ en un drapeau complet $\wp$ tel que $f$
se rel{\`e}ve en $\ft:C\to\Xt_\wp$ avec $\ft_*[C]\in\app$.
\end{prop}

\begin{preu}
 Remarquons tout d'abord que pour qu'un morphisme $f:C\to X$ se rel{\`e}ve
dans $\Xt_\wp$, il suffit qu'il rencontre l'ouvert $U_0$. En effet,
l'image d'un ouvert de $C$ rencontre alors $U_0$ et l'application
$\ft$ est d{\'e}finie au moins sur cet ouvert. Comme $C$ est une courbe
lisse, il se prolonge {\`a} $C$ toute enti{\`e}re.

Par ailleurs remarquons {\'e}galement qu'un {\'e}l{\'e}ment $x\in U$ est dans
$U_0$ si et seulement si son image $p(x)$ n'est contenue dans aucun
sous-espace du drapeau complet de $V$ obtenu comme la trace de $\wp$. 

Remarquons enfin que le choix d'un drapeau complet dans $V$ est
{\'e}quivalent au choix d'un drapeaux complet de l'espace $\C^{2n}$
completant $V$. Nous montrons le lemme suivant :

\begin{lemm}
Soit $g:C\to\p(V)$ un morphisme vers $\p(V)$ dont l'image r{\'e}duite
n'est contenue dans  aucun sous-espace lin{\'e}aire de codimension 2,
alors on peut compl{\'e}ter $V$ en un drapeau complet $\wp$ tel que $g$ se
rel{\`e}ve en
$\gt:C\to\pt_{\wp}(V)$ avec $\gt_*[C]\cdot\xi_i>0 \textrm{ pour tout }
i\in[1,n-2]$ et $\gt_*[C]\cdot\xi_{n-1}\geq0$.
\end{lemm}

\begin{preu}
  Si $0=F_0\subset F_1\subset\cdots\subset F_{n-1}\subset V$ est le
drapeau complet de $V$ induit par $\wp$, alors la vari{\'e}t{\'e} $\pt_\wp(V)$
est d{\'e}crite par 
$$\pt_\wp(V)=\left\{(V_{n-i})_{i\in[1,n-1]}\in\prod_{i=1}^{n-1}\G(n-i,V)\ /\
F_{i-1}\subset V_i\subset V_{i+1} \textrm{ pour tout }
i\in[1,n-1]\right\}$$ 
avec $V_{n}=V$. La projection $\pi':\pt_\wp(V)\to\p(V)$ est donn{\'e}e par
la projection sur le dernier facteur $(V_i)_{i\in[1,n-1]}\mapsto
V_1$. Chaque diviseur $\xi_i$ est donn{\'e} par l'{\'e}quation
$V_{n-i}=F_{n-i}$. 

On construit par r{\'e}currence une suite de points $(P_i)_{i\in[1,n-1]}$
de $C$ et les sous-espaces $F_i$ du drapeau complet de telle sorte que
pour tout $i\in[1,n-1]$ on ait :
$$\left\{
  \begin{array}{l}
g(P_i)\in F_{i+1}\\
g(P_i)\not\in F_i.\\
  \end{array}\right.$$

Si une telle donn{\'e}e existe, alors le point $g(P_{n-1})$ est dans l'ouvert
o{\`u} $\pi'$ est bijective donc on peut relever $g$ en $\gt$ avec par
ailleurs 
$$\gt(P_{n-1})=(g(P_{n-1})+F_{n-2},\cdots,g(P_{n-1})+F_1,g(P_{n-1})).$$
En particulier, le point $\gt(P_{n-1})$ n'est contenu dans aucun des
$\xi_i$ et on a donc $\gt_*[C]\cdot\xi_i\geq0$ pour tout
$i\in[1,n-1]$. 
Par ailleurs, on voit que le point
$\gt(P_i)=(V_{n-j}(P_j))_{j\in[1,n-1]}$ v{\'e}rifie
$V_{i+1}(P_i)=g(P_i)+F_i=F_{i+1}$. Le point $\gt(P_i)$ est donc dans le
diviseur $\xi_{i+1}$. On a l'intersection souhait{\'e}e.

Il reste donc {\`a} construire les points $P_i$ et les sous-espaces
$F_i$. On prend les points $P_i$ en position g{\'e}n{\'e}rale sur $C$. En
particulier, l'espace engendr{\'e} par $g(P_1),\cdots,g(P_i)$ doit {\^e}tre
de dimension 
$i$ pour tout $i\in[1,n-1]$. On construit les $F_i$ par r{\'e}currence
descendante sur $i$. Le sous-espace de dimension $n$ du drapeau est
$V$. Une fois $F_{i+1}$
fix{\'e} contenant $g(P_1),\cdots, g(P_i)$, on choisit $F_{i}$ contenu dans
$F_{i+1}$ et contenant $g(P_1),\cdots, g(P_{i-1})$ mais ne contenant pas
$g(P_i)$. Ceci est possible par hypoth{\`e}se sur les $P_i$ et on a une
famille de dimension 1 de choix pour $F_i$.
Ainsi la donn{\'e}e des $P_i$ et des $F_i$ existe et il en existe une
famille de dimension $2(n-1)$. Il existe donc une famille de dimension
$2n-3$ de drapeaux qui conviennent (une fois le drapeau fix{\'e}, il y a
un nombre fini de choix sur les points $P_i$ sauf pour $x_{n-1}$). 
\end{preu}

Le lemme nous permet de fixer un drapeau complet $\wp$ tel que $g=p\circ
f$ se rel{\`e}ve en $\gt:C\to\pt_\wp(V)$ avec $\gt_*[C]\cdot\xi_i>0$ pour
tout $i\in[1,n-2]$. 

Par ailleurs, comme $p\circ f(P_{n-1})$ n'est contenu dans aucun des
sous-espace du drapeau induit par $\wp$ dans $V$, le point
$f(P_{n-1})$ est dans l'ouvert $U_0$. La fl{\`e}che $f$ se rel{\`e}ve donc en
$\ft:C\to\Xt_\wp$ et on a $\gt=q\circ\ft$. Le fait que $f$ rencontre
$U_0$, impose que pour tout $i\in[1,r]$, on a
$\ft_*[C]\cdot\xi_i\geq0$. Les intersections $\gt_*[C]\cdot\xi_i>0$ pour
tout $i\in[1,n-2]$ imposent que l'on a $\ft_*[C]\cdot\xi_i>0$ pour
tout $i\in[1,n-2]$.
\end{preu}

\section{Irr{\'e}ductibilit{\'e}}

Nous donnons dans ce paragraphe la preuve de notre r{\'e}sultat
principal.

\subsection{Courbes g{\'e}n{\'e}rales}

Nous commen{\c c}ons par montrer que les morphismes qui ne satisfont pas
les conditions de la proposition \ref{intersection} ne peuvent former
une composante irr{\'e}ductible du sch{\'e}ma des morphismes pour un degr{\'e}
assez grand. Rappelons que le morphisme $p:U\to\p(V)$ est le morphisme
de projection du fibr{\'e} vectoriel $\Lambda^2T_{\p(V)}(-1)$.

\begin{prop}
\label{reste}
  Soit $\a\in A_1(X)$ une classe de degr{\'e} $d\geq n-1$. Le ferm{\'e} du
  sch{\'e}ma des morphismes $\MorC{i^*\a}{U}$ form{\'e} des fl{\`e}ches $f$ telles
  que la courbe $p\circ f(C)$ est contenue dans un sous-espace
  lin{\'e}aire de codimension 2 
ne peut former une
composante irr{\'e}ductible de $\MorC{i^*\a}{U}$.
\end{prop}

\begin{preu}
  Nous consid{\'e}rons le morphisme
$\MorC{i^*\a}{U}\to\MorC{{p}_*i^*\a}{\p(V)}$ induit par
$p$. Pour d{\'e}finir ${p}_*$ nous consid{\'e}rons les classes de
1-cycles comme des {\'e}l{\'e}ments du dual du groupe de Picard et utilisons
la transpos{\'e}e de ${p}^*$. L'image du ferm{\'e} consid{\'e}r{\'e} dans
$\MorC{{p}_*i^*\a}{\p(V)}$ est stratifi{\'e} par les familles
$(M_k)_{k\in[2,n-1]}$ des morphismes $f:C\to\p(V)$ dont l'image
est contenue dans un sous-espace lin{\'e}aire de codimension
$k$. La dimension de $M_k$ est $(n-k)d+k(n-k)$.

La fibre du morphisme $\MorC{i^*\a}{U}\to\MorC{{p}_*i^*\a}{\p(V)}$ au
dessus de $f:C\to\p(V)$ est donn{\'e}e par
$H^0(f^*\Lambda^2T_{\p(V)}(-1))$. Pour un morphisme dont l'image est
contenue dans un sous-espace lin{\'e}aire de codimension $k$ mais pas dans
un sous-espace lin{\'e}aire de codimension plus grande, on a 
$$f^*T_{\p(V)}(-1)=\oo_C^k\oplus E$$
o{\`u} $E$ est un fibr{\'e} de rang $n-1-k$, de degr{\'e} $d$ sur $C$, engendr{\'e}
par ses sections (car $T_{\p(V)}(-1)$ l'est) et sans facteur trivial. 
Ainsi le faisceau $f^*\Lambda^2T_{\p(V)}(-1)$ est isomorphe {\`a}
la somme directe suivante:
$$\oo_C^{\frac{k(k-1)}{2}}\oplus E^k\oplus \Lambda^2E.$$ 
Mais $E$ est engendr{\'e} par ses sections et on a $H^1E=0$ ce qui
impose que c'est {\'e}galement le cas de $\Lambda^2E$ (on a une
surjection de $H^0E\otimes E$ vers $\Lambda^2E$). 
Le groupe $H^1(C,f^*\Lambda^2T_{\p(V)}(-1))$ est donc de
dimension $\displaystyle{\frac{k(k-1)}{2}}$ et le groupe
  $H^0(C,f^*\Lambda^2T_{\p(V)}(-1))$ est de dimension
$$(n-2)d+\displaystyle{\frac{k(k-1)}{2}}.$$

La famille des morphismes de $\MorC{i^*\a}{U}$ qui s'envoie dans $M_k$
est donc de dimension 
$$(2n-k-2)d+k(n-k)+\frac{k(k-1)}{2}.$$

Enfin la dimension d'une composante irr{\'e}ductible de $\MorC{i^*\a}{U}$
est au moins $2(n-1)d$ et on a
$$2(n-1)d>(2n-k-2)d+k(n-k)+\frac{k(k-1)}{2}$$
d{\`e}s que $d\geq n-1$ et $k\geq2$.
\end{preu}

\subsection{L'incidence}

Soit $B$ un sous-groupe de Borel de $G$, la vari{\'e}t{\'e} des drapeaux
complets est $G/B$. 
Nous avons vu qu'{\`a} chaque drapeau complet $\wp$ correspond une
r{\'e}solution de Bott-Samelson $\Xt_\wp\to X$. Consid{\'e}rons l'incidence
suivante :
$$
\xymatrix{
I  \ar@{->}[r]^-{\pr_2} \ar[d]_{\pr_1} & {\MorC{\a}{X}}\\
G/B &
\\}
$$
o{\`u} $I$ est form{\'e}e des couples $(\wp,f)\in G/B\times\MorC{\a}{X}$ tels
que $f(C)$ est contenue dans l'ouvert $U$ associ{\'e} au drapeau $\wp$, la
fl{\`e}che $f$ se rel{\`e}ve en $\ft:C\to\Xt_\wp$ et on a
$\ft_*[C]\in\app$ (voir paragraphe \ref{irredu} pour la d{\'e}finition
de cet ensemble). 

Les propositions \ref{intersection} et \ref{reste} montrent que la
fl{\`e}che $\pr_2:I\to\MorC{\a}{X}$ est dominante d{\`e}s que $d\geq n-1$. Rappelons
{\'e}galement que nous avons vu dans la preuve de la proposition
\ref{intersection} que la dimension de la fibre de $\pr_2$ au-dessus
d'un point g{\'e}n{\'e}ral est $2(n-1)+\displaystyle{\frac{n(n-1)}{2}}$ (le
second terme provenant du choix de $V$).

\subsection{L'irr{\'e}ductibilit{\'e}}

Par d{\'e}finition de $I$, on a une fl{\`e}che surjective :
$$\coprod_{\at\in\app,\ \pi_*(\at)=\a}G/B\times\MorC{\at}{\Xt_\wp}\to
I$$
et par composition une fl{\`e}che dominante
$$\coprod_{\at\in\app,\ \pi_*(\at)=\a}G/B\times\MorC{\at}{\Xt_\wp}\to
\MorC{\a}{X}.$$

On a vu {\`a} la section \ref{dimension} que pour tout {\'e}l{\'e}ment
$\at\in\app$, la dimension de $\MorC{\at}{\Xt_\wp}$ est inf{\'e}rieure ou
{\'e}gale {\`a}
$$2(n-1)d-\frac{(n-2)(n-3)}{2}$$
avec {\'e}galit{\'e} si et seulement si on a $\at=\at_0$. Pour une
classe $\at\in\app$ fix{\'e}e telle que $\pi_*\at=\a$, l'image de
$G/B\times\MorC{\at}{\Xt_\wp}$ dans $\MorC{\a}{X}$ est donc de
dimension inf{\'e}rieure {\`a} 
$$2(n-1)d-\frac{(n-2)(n-3)}{2}+\dim(G/B)-\dim(\pr_2^{-1}(f))\leq
2(n-1)d.$$
Cependant une composante irr{\'e}ductible de $\MorC{\a}{X}$ est de
dimension au moins $2(n-1)d$, ainsi pour $\at\neq\at_0$ l'image de
$G/B\times\MorC{\at}{\Xt_\wp}$ dans $\MorC{\a}{X}$ ne peut contenir de
composante irr{\'e}ductible. On a donc une application dominante 
$$G/B\times\MorC{\at_0}{\Xt_\wp}\to \MorC{\a}{X}.$$
Le sch{\'e}ma $\MorC{\at_0}{\Xt_\wp}$ {\'e}tant irr{\'e}ductible, on en d{\'e}duit le
r{\'e}sultat principal de cet article :

\begin{theo}
  Pour $\a\in A_1(X)$ tel que $d=\deg(\a)\geq n-1$, le sch{\'e}ma
  $\MorC{\a}{X}$ est irr{\'e}ductible de dimension $2(n-1)d$.
\end{theo}

\section{Appendice}

Nous tra{\c c}ons ici le carquois associ{\'e} {\`a} l'{\'e}criture r{\'e}duite donn{\'e}e au
paragraphe \ref{Bott}. Ce carquois d{\'e}pend de la parit{\'e} de $n$. La
fl{\`e}che $\b$ de l'ensemble des sommets vers les racines simples est ici
donn{\'e}e par la projection verticale. Nous n'avons pas trac{\'e} les fl{\`e}ches
sur les ar{\^e}tes pour ne pas surcharger le dessin, mais toutes les
fl{\`e}ches vont vers le bas.

\psset{xunit=1cm}
\psset{yunit=1cm}
\centerline{\begin{pspicture*}(-0.2,-8.5)(7.5,7)
\put(3,-8){$n$ pair.}
\put(7,6.2){1}
\put(5,5.2){2}
\put(0,0.2){$n-1$}
\put(6,4.2){$n$}
\put(5,3.2){$n+1$}
\put(1,-0.8){$2n-3$}
\put(7,-5.8){$r$}
\psline(0,0)(1.5,1.5)
\psline[linestyle=dotted](1.5,1.5)(3.5,3.5)
\psline(3.5,3.5)(5,5)
\psline(5,5)(7,6)
\psline(1,1)(2.5,-0.5)
\psline(5,5)(6,4)
\psline(4,4)(5,3)
\psline(5,3)(7,2)
\psline(6,4)(4.5,2.5)
\psline(7,2)(6,1.5)
\psline[linestyle=dotted](6,1.5)(5,1)
\psline[linestyle=dotted](4.5,2.5)(2.5,0.5)
\psline[linestyle=dotted](3,3)(5,1)
\psline(7,6)(5,3)
\psline(6,4)(5.5,2.5)
\psline[linestyle=dotted](5,1)(5.5,2.5)
\psline(7,-2)(5,-5)
\psline[linestyle=dotted](5,0.7)(5,-0.7)
\psline(0,0)(1.5,-1.5)
\psline[linestyle=dotted](1.5,-1.5)(3.5,-3.5)
\psline(3.5,-3.5)(5,-5)
\psline(5,-5)(7,-6)
\psline(1,-1)(2.5,0.5)
\psline(5,-5)(6,-4)
\psline(4,-4)(5,-3)
\psline(5,-3)(7,-2)
\psline(6,-4)(4.5,-2.5)
\psline(7,-2)(6,-1.5)
\psline[linestyle=dotted](6,-1.5)(5,-1)
\psline[linestyle=dotted](4.5,-2.5)(2.5,-0.5)
\psline[linestyle=dotted](3,-3)(5,-1)
\put(-0.1,-0.1){$\bullet$}
\put(0.9,0.9){$\bullet$}
\put(3.9,3.9){$\bullet$}
\put(4.9,4.9){$\bullet$}
\put(6.9,5.9){$\bullet$}
\put(5.9,3.9){$\bullet$}
\put(4.9,2.9){$\bullet$}
\put(1.9,-0.1){$\bullet$}
\put(0.9,-1.1){$\bullet$}
\put(6.9,1.9){$\bullet$}
\put(3.9,-4.1){$\bullet$}
\put(4.9,-5.1){$\bullet$}
\put(6.9,-6.1){$\bullet$}
\put(5.9,-4.1){$\bullet$}
\put(6.9,-2.1){$\bullet$}
\put(4.9,-3.1){$\bullet$}
\put(-0.1,-7.1){$\bullet$}
\put(0.9,-7.1){$\bullet$}
\put(1.9,-7.1){$\bullet$}
\put(6.9,-8.1){$\bullet$}
\put(5.9,-6.1){$\bullet$}
\put(4.9,-7.1){$\bullet$}
\put(3.9,-7.1){$\bullet$}
\psline{->}(3,-3.5)(3,-6)
\psline(0,-7)(2.5,-7)
\psline[linestyle=dotted](2.5,-7)(3.5,-7)
\psline(3.5,-7)(5,-7)
\psline(6,-6)(5,-7)
\psline(7,-8)(5,-7)
\put(2.5,-4.8){$\b$}
\end{pspicture*}
\begin{pspicture*}(-0.2,-8.5)(7.5,7)
\put(3,-8){$n$ impair.}
\put(6,4.2){$1$}
\put(5,3.2){$2$}
\put(1,-0.8){$n-1$}
\put(7,-5.8){$r$}
\psline(2,0)(2.5,-0.5)
\psline(5,3)(7,2)
\psline(6,4)(4.5,2.5)
\psline(7,2)(6,1.5)
\psline[linestyle=dotted](6,1.5)(5,1)
\psline[linestyle=dotted](4.5,2.5)(2.5,0.5)
\psline[linestyle=dotted](4,2)(5,1)
\psline(6,4)(5.5,2.5)
\psline[linestyle=dotted](5,1)(5.5,2.5)
\psline(7,-2)(5,-5)
\psline[linestyle=dotted](5,0.7)(5,-0.7)
\psline(1,-1)(1.5,-1.5)
\psline[linestyle=dotted](1.5,-1.5)(3.5,-3.5)
\psline(3.5,-3.5)(5,-5)
\psline(5,-5)(7,-6)
\psline(1,-1)(2.5,0.5)
\psline(5,-5)(6,-4)
\psline(4,-4)(5,-3)
\psline(5,-3)(7,-2)
\psline(6,-4)(4.5,-2.5)
\psline(7,-2)(6,-1.5)
\psline[linestyle=dotted](6,-1.5)(5,-1)
\psline[linestyle=dotted](4.5,-2.5)(2.5,-0.5)
\psline[linestyle=dotted](3,-3)(5,-1)
\psline{->}(3,-3.5)(3,-6)
\psline(1,-7)(2.5,-7)
\psline[linestyle=dotted](2.5,-7)(3.5,-7)
\psline(3.5,-7)(5,-7)
\psline(6,-6)(5,-7)
\psline(7,-8)(5,-7)
\put(5.9,3.9){$\bullet$}
\put(4.9,2.9){$\bullet$}
\put(1.9,-0.1){$\bullet$}
\put(0.9,-1.1){$\bullet$}
\put(6.9,1.9){$\bullet$}
\put(0.9,-7.1){$\bullet$}
\put(1.9,-7.1){$\bullet$}
\put(6.9,-8.1){$\bullet$}
\put(5.9,-6.1){$\bullet$}
\put(4.9,-7.1){$\bullet$}
\put(3.9,-7.1){$\bullet$}
\put(3.9,-4.1){$\bullet$}
\put(4.9,-5.1){$\bullet$}
\put(6.9,-6.1){$\bullet$}
\put(5.9,-4.1){$\bullet$}
\put(6.9,-2.1){$\bullet$}
\put(4.9,-3.1){$\bullet$}
\put(2.5,-4.8){$\b$}
\end{pspicture*}}

\begin{small}

\vs 0.2 cm

\noi
{\textsc{Universit{\'e} Pierre et Marie Curie - Paris 6}}

\vs -0.1 cm

\noi
{\textsc{UMR 7586 --- Institut de Math{\'e}matiques de Jussieu}}

\vs -0.1 cm

\noi
{\textsc{175 rue du Chevaleret}}

\vs -0.1 cm

\noi
{\textsc{75013 Paris,}} \hs 0.2 cm{\textsc{France.}}

\vs -0.1 cm

\noi
{email : \texttt{nperrin@math.jussieu.fr}}

\end{small}

\end{document}